\documentclass[a4paper,12pt]{article}
\usepackage[latin1]{inputenc}
\usepackage[T1]{fontenc}
\usepackage[english]{babel}
\usepackage{amssymb}
\usepackage{latexsym}
\usepackage{xcolor}
\usepackage{amsmath}
\usepackage{mathtools}
\usepackage{amsthm}
\usepackage{mathrsfs}
\usepackage{amsxtra}
\usepackage{eucal}
\usepackage[left=1.5cm,right=1.5cm,top=1cm,bottom=1.5cm]{geometry}
\usepackage{tikz}
\usetikzlibrary{patterns}
\usepackage{amsfonts}
\usepackage{lmodern}
\usepackage{graphicx}
\usepackage{float}
\usepackage{hyperref}
\providecommand{\LyX}{L\kern-.1667em\lower.25em\hbox{Y}\kern-.125emX\@}

\def\R{\mathbb{R}}
\def\E{\mathbb{E}}
\def\dint{\displaystyle\int}

\newtheorem{thm}{Theorem}[section]

\newtheorem{cor}[thm]{Corollary}

\newtheorem{examples}[thm]{Examples}

\newtheorem{lem}[thm]{Lemma}

\newtheorem{prop}[thm]{Proposition}
\newtheorem{rk}[thm]{Remark}

\newenvironment{class}[1][AMS Classification]{\textbf{#1.} }{}
\newenvironment{MC}[1][Key words]{\textbf{#1.} }{}

\def\al{\alpha}
\def\si{\sigma}
\def\be{\beta}

\title{\textbf{ On  Carath\'eodory approximate scheme for a class of one-dimensional doubly perturbed diffusion processes}}

\author{ R. BELFADLI  $^{1}$
$\quad$ L. BOULANBA   $^{2}\quad$ Y. OUKNINE $^{3}$ $^{4}\quad$
\vspace*{0.1in}\\
$^{1}$  Department of Mathematics, Faculty
of Sciences and Technologies,\\ Cadi Ayyad University 2390 Marrakesh, Morocco.\\
rachid.belfadli@uca.ac.ma 
\vspace*{0.1in}\\
$^{2}$ Department of Mathematics, CRMEF SM, Morocco.\\ Laboratoire LISTI, ENSA Agadir, Universit\'e Ibn Zohr.\\
\quad l.boulanba@crmefsm.ac.ma 
\vspace*{0.1in}\\
$^{3}$ Cadi Ayyad University, Faculty of Sciences Semlalia\\
 Department of Mathematics, B.P. 2390, Marrakesh, 40.000, Morocco.\\
ouknine@uca.ac.ma
\vspace*{0.1in}\\
$^{4}$ Mohammed VI Polytechnic University, Africa Business School, \\Avenue Mohammed Ben Abdellah Regragui,\\
Madinat  Al Irfane, BP 6380, RABAT, Morocco. \\
 youssef.ouknine@um6p.ma
}

\begin{document}
\maketitle \textbf{Abstract.} In this paper, we introduce and study the convergence of new Carath\'eodory's approximate solution for one-dimensional   $\alpha, \beta$-doubly perturbed stochastic differential equations (DPSDEs) with parameters  $\alpha <1$ and $\beta <1$   such that $|\rho|  < 1$, where $ \rho : = \frac{\alpha\beta}{(1-\alpha)(1-\beta)}$. Under Lipschitz's condition on the coefficients, we establish the  $L^{p}$-convergence of the  Carath\'{e}odory approximate solution uniformly in time, for all $p\geq 2$. As a consequence, and relying only on our scheme,  we obtain the existence and uniqueness of strong solution for $\alpha, \beta$-DPSDEs. Furthermore, an extension to non-Lipschitz coefficients are also studied. Our results  improve  earlier work by Mao and al. \cite{MHM2018}.
\\

 \begin{MC}
Doubly perturbed stochastic differential equations; Carath\'eodory numerical scheme. Skorohod's lemma.
\end{MC}

 \begin{class}  Primary 60H10, 60H35; Secondary 60J60.
\end{class}
\section{Introduction}
 
Let $T>0$ be a fixed time horizon, and  consider the following one-dimensional  $\alpha, \beta$-doubly perturbed stochastic differential equations (SDEs): 
\begin{eqnarray}\label{eqbb1}
X_t=x_{0} + \int_0^t b(s, X_s)ds + \int_0^t \si(s, X_s)dW_s  +\al\max_{0\leq
s\leq t}X_s +\be \min_{0\leq s \leq t}X_s, \, \, t \in [0, T],
\end{eqnarray}
where $x_0 \in \mathbb{R}$, $W=\left(W_{t}\right)_{t\geq 0}$ is a standard Brownian motion defined on a filtered probability space $(\Omega, \mathcal{F}, (\mathcal{F}_t)_{t\geq 0}, \mathbb{P})$ with a filtration satisfying the usual conditions.

 Stochastic differential equations  involving the past maximum and/or the past minimum processes has attracted  much interest and has been  studied  by several authors during the past few decades, (see for instance \cite{LGY1990, CPY95, Davis1996, PW1997, CD1999, CDH2000, Doney-Zhang-2005, Belfadli2009, Belfadli2024, MHM2018, GM2023}).

Assuming Lipschitz's conditions  for the real valued functions $b$ and $\sigma$, existence and uniqueness  of a strong solution $X=(X_t)_{t\geq 0}$  is provided in  \cite{CPY95, Davis1996, Belfadli2009} whenever the real parameters  $\alpha$ and $\beta$ are such that:
\begin{equation}\label{hyp}
 \alpha <1, \,\,   \,\, \beta <1 \, \, \,\mbox{and}\,\,\,  |\rho|:= \frac{|\alpha \beta|}{(1- \alpha)(1-\beta)}|  < 1.
\end{equation}

In their paper \cite{MHM2018}, Mao \textit{et al}. have introduced and studied Carath\'eodory approximation scheme  for the SDE (\ref{eqbb1}). Actually, unlike the Picard iterative method, it is  an approximation  that avoid the need of   calculating at each step the previous  values of the sequence. Namely, the authors introduce the following  approximate solutions $X^{n}: [-1, T] \longrightarrow \mathbb{R}$ defined, for all integer $n \geq 1$, by
\begin{equation}\label{eq3:appro}
\left\{
\begin{array}{l}
 X_t^{n} = x_0,  \,\,\, \, \, \mbox{for}\,\,\,  -1 \leq t \leq 0, \\
X_t^{n} :=x_{0} + \int_0^t b(s, X^{n}_{\eta_n(s)})ds + \int_0^t \si(s, X^{n}_{\eta_n(s)})dW_s \\
\,\,\,\,\qquad\qquad +\al\max_{0\leq
s\leq t}X^{n}_{\eta_n(s)} +\be \min_{0\leq s \leq t}X^{n}_{\eta_n(s)}, \, \,\, \, \mbox{for}\,\,  t \in (0, T],
\end{array}
\right.
\end{equation}
where $\eta_n(s):=s-\frac{1}{n}$ for $s\geq 0$.   Notice that $X^{n}_t$ can be calculated step by step on each interval $[0, \frac{1}{n})$, $[\frac{1}{n}, \frac{2}{n}), \ldots$

Still in \cite{MHM2018}, it is proved that the scheme $\{X^{n}\}_{n\in \mathbb{N}^{\ast}}$ verifies, under Lipschitz conditions on $b$ and $\sigma$ and when  $\sup_{0\leq t \leq T}(|b(t, 0)| +|\sigma(t, 0)|)< \infty $, that 
\begin{equation}\label{eq:1}
\E[\sup_{0\leq t \leq T}|X^{n}_t-X_t|^2] \leq \frac{c}{n},
\end{equation}

where $c$ is a constant independent of $n$. In addition,   when $b$ and $\sigma$ are no more Lipschitz and under appropriate conditions, it is proved that  the following convergence holds
\begin{equation}\label{eq:2}
\lim_{n\rightarrow +\infty} \E[\sup_{0\leq t \leq T}|X^{n}_t-X_t|^2] =0.
\end{equation}

Unfortunately, the proof of these two results requires  the parameters $\alpha$ and $\beta$ to satisfy  $|\alpha| + |\beta| <1$. Of course, it should come as no surprise that this assumption is quite restrictive if compared with the aforementioned condition  leading to the existence and uniqueness result.  Indeed, a careful examination of the arguments given in the proof of  (\ref{eq:1}) and (\ref{eq:2})  show that this is due heavily to the way the Carath\'eodory scheme (\ref{eq3:appro}) is defined. Thus, it is sensible to modify the scheme (\ref{eq3:appro}) in order to improve the above results and go beyond the  condition $|\alpha| + |\beta|<1$.

Our principal  motivation in writing this paper 
is to introduce  a more appropriate Carath\'eodory scheme  which will be very convenient for establishing the above results for  all parameters $\alpha$ and $\beta$ satisfying  (\ref{hyp}).

We first give an explanation of the key idea behind our new approximation.  Using Skorohod's lemma ( cf. \cite{Revuz-Yor-2005}, p. $239$ or \cite{KS1988}, p. $210$) to the equation  (\ref{eqbb1}), when $x_0=0$, yields :
\begin{equation}\label{eq: 3}
M_t^X : = \frac{1}{1 - \alpha}\max_{0 \leq s \leq t}\left\{ \left(  \int_{0}^{s}\sigma\left(u, X_{u}\right)dW_u +  \int_{0}^{s}b\left(u, X_{u}\right)du + \beta I^X_{s} \right)^{+}\right\}, 
\end{equation}
and 
\begin{equation}\label{eq: 4}
I_t^X : = \frac{1}{\beta- 1}\max_{0 \leq s \leq t}\left\{ \left( \int_{0}^{s}\sigma\left(u, X_{u}\right)dW_u - \int_{0}^{s}b\left(u, X_{u}\right)du -\alpha M^X_{s} \right)^{+}\right\} 
\end{equation}
where $M^{X}_t:= \max_{0\leq s \leq t}X_s$ and  $I^{X}_t:= \min_{0\leq s \leq t}X_s$ and  $x^{+} = \max(x, 0)$.

 The nice thing about the identities (\ref{eq: 3}) and (\ref{eq: 4}) that we would like to highlight is that one can use the whole past of the process $\{I^{X}_s, 0\leq s \leq t\}$ to evaluate $M^{X}_t$ as well as the past of  $\{M^{X}_s, 0\leq s \leq t\}$ to calculate $I^{X}_t$. Taking into account this observation, it seems natural to introduce the  following $\mathbb{R}^{3}$-valued Carath\'eodory scheme  $(X_{t}^{n}, M^{n}_t, I^{n}_t)_{n \geq 1}$  associated to (\ref{eqbb1}) with $x_0=0$. For each integer $n\geq 1$, 
\begin{equation*}
 X_t^n =  M_t^n = I_t^n = 0,\,\,\,\, \text{for} \ -1 \leq t \leq 0,
\end{equation*}
and
\begin{equation} \label{eqbb1bis}
 X_t^n  =   \int_{0}^{t}\sigma\left(s, X^n_{\eta_n(s)}\right)dW_s +  \int_{0}^{t}b\left(s, X^n_{\eta_n(s)}\right)ds + \alpha M_t^n  + \beta I_t^n, \text{for}\\,\,  t\in (0, T],
\end{equation}
where $M_t^n$  and $I_t^n$ are respectively given by 
\begin{equation}\label{eqbb3}
M_t^n : = \frac{1}{1 - \alpha}\max_{0 \leq s \leq t}\left\{\left(\int_{0}^{s}\sigma\left(u, X^n_{\eta_n(u)}\right)dW_u +  \int_{0}^{s}b\left(u, X^n_{\eta_n(u)}\right)du + \beta I^n_{\eta^{+}_n(s)} \right)^{+}\right\}, 
\end{equation}
\begin{equation}\label{eqbb4}
I_t^n : = \frac{1}{\beta -1}\max_{0 \leq s \leq t}\left\{\left(- \int_{0}^{s}\sigma\left(u, X^n_{\eta_n(u)}\right)dW_u -  \int_{0}^{s}b\left(u, X^n_{\eta_n(u)}\right)du -\alpha M^n_{\eta^{+}_n(s)}\right)^{+}\right\}.
\end{equation}

\noindent Our choice of $M^{n}$ and $I^{n}$ is heuristically explained by (\ref{eq: 3}) and (\ref{eq: 4}). Obviously, because of  (\ref{eqbb3}) and (\ref{eqbb4}),  $X^{n}_t$ is defined 
repeatedly step by step on each interval $[-1, 0]$,  $[0, \frac{1}{n}[$, $[\frac{1}{n}, \frac{2}{n}[$,  $\cdots$, without using any of the previous approximations $X^{i}, i =1, \cdots, n-1$. That is, $(X^n)$ is a Carath\'{e}odory approximation for (\ref{eqbb1}) which is continuous.

 To simplify notation, we set 
 \begin{equation}\label{eq:Phin}
 \displaystyle  \Phi_t^n : = \int_{0}^{t}\sigma\left(s, X^n_{\eta_n(s)}\right)dW_s + \int_{0}^{t}b\left(s, X^n_{\eta_n(s)}\right)ds.
\end{equation}  
  Inserting (\ref{eqbb4}) into (\ref{eqbb3}) we obtain
\begin{equation}\label{eqbb7}
M_t^n = \frac{1}{1 - \alpha}\max_{0 \leq s \leq t}\left\{\left(\Phi^n_{s} +\frac{\beta}{\beta - 1}\max_{0\leq u \leq {\eta^{+}_n(s)}}\left(-\Phi^n_{u} -\alpha M^n_{\eta^{+}_n(u)}\right)^{+}\right)^{+}\right\}, 
\end{equation}  
and similarly
\begin{equation}\label{eqbb7bis}
I_t^n = \frac{1}{\beta - 1}\max_{0 \leq s \leq t}\left\{\left(-\Phi^n_{s} -\frac{\alpha}{ 1-\alpha}\max_{0\leq u \leq {\eta^{+}_n(s)}}\left(\Phi^n_{u} +\beta I^n_{\eta^{+}_n(u)}\right)^{+}\right)^{+}\right\}.
\end{equation}  

\noindent Remark that unlike expressions  (\ref{eqbb3}) and (\ref{eqbb4}),  the formulae (\ref{eqbb7}) and  (\ref{eqbb7bis}) feature the processes $M^{n}$ and $I^{n}$ both on   their left and right hand sides.  It is worth emphasizing  that the representations (\ref{eqbb7}) and (\ref{eqbb7bis})  will play  prominent role in our study. Indeed, as it will be shown later, most estimates  related to $M_t^n$ and $I_t^n$ are obtained by using  (\ref{eqbb7}) and (\ref{eqbb7bis}).

The remainder of this paper is structured as follows.  In Section 2, we establish the existence  and uniqueness result of solution to  (\ref{eqbb1}) and show that the Carath\'{e}odory approximate solution converges to the solution of equation (\ref{eqbb1}) under the global Lipschitz condition. While in Section 3, we extend the existence and convergence results of Section 2 to the case of equation (\ref{eqbb1}) with non-Lipschitz coefficients.

Throughout this paper, we assume the parameters $\alpha$ and $\beta$ satisfy the condition (\ref{hyp}). We will also use the letter $c>0$ to indicate a constant whose value is immaterial that may vary from line to line and may depend only on $\alpha$, $\beta$, $T$ and $p$.

\section{Statement of the main result}

Here we show,  under Lipschitz condition  on $b$ and $\sigma$, the convergence of the approximation $(X^n)$ to the unique solution of  (\ref{eqbb1}). More precisely, let us consider the following assumption: 
\vspace{.4cm}

\noindent \textbf{($\mathbf{H}1$)}: there exists a constant $K>0$, such that
\begin{equation*}
\label{eq2} \left\{
\begin{array}{l}
|\sigma(t,x) - \sigma(t,y)|  + |b(t,x) - b(t,y)|  \leq K |x-y|,\,\,\, \\
|\sigma(t,0)|  + |b(t, 0)|\leq K, 
\end{array}
\right.
\end{equation*}
for all $ t \in [0, T]$, and   $x,  y  \in \mathbb{R}.$
\noindent Under the above condition  \textbf{($\mathbf{H}1$)}, we see that for any $x\in  \mathbb{R}$ and $p\geq 1$,
\begin{equation}\label{estimate:1}
\sup_{t\in [0, T]} |\sigma(t,x)|^{p}  \leq c ( 1 + |x|^p)\, \, \, \mbox{and}\,\, \, \sup_{t\in [0, T]} |b(t,x)|^{p} \leq c ( 1 + |x|^p).
\end{equation}

\vspace{0.2cm}

\noindent Now, we state the first main result of this paper.
\begin{thm}\label{thmbb1} Assume that the functions $b$ and $\sigma$ are continuous and  satisfy  \textbf{($\mathbf{H}1$)}. Then the sequence $\{X^n\}$ defined by (\ref{eqbb1bis}) converges in $\mathrm{L}^{p}(\Omega, C([0,T];\mathbb{R}))$ and almost surely uniformly on $[0, T]$. The limit is the unique solution to the SDE (\ref{eqbb1}) with $x_0=0$. Furthermore, there exists a constant $c>0$ such that
\begin{equation} \label{rate:1}
\mathbb{E}\left(\sup_{t\in [0, T]}|X^n_t - X_t|^{p}\right) \leq \frac{c}{n^{p\gamma}},
\end{equation}
for any $p>2$ and $\gamma \in (0, \frac{1}{2}-\frac{1}{p})$.
\end{thm}
\begin{rk}
It is worth nothing that, based on the scheme (\ref{eqbb1bis}), (\ref{eqbb3}) and (\ref{eqbb4}), Theorem \ref{thmbb1} is a sharpened version of Theorem $2. 1$ in Mao and al. \cite{MHM2018} for all $p> 2$ and parameters $\alpha <1$ and $\beta <1$ satisfying  only $|\rho|<1$ rather than $|\alpha| + |\beta|<1$.
\end{rk}
 In order to prove the first part of Theorem \ref{thmbb1}, it suffices to show that  $(X^n)$ is a Cauchy sequence in $\mathrm{L}^{p}(\Omega, C([0,T];\mathbb{R}))$. As a preliminary, we will need three auxiliary lemmas. 
\vspace{0.2cm}
\begin{lem} \label{lembb1}
Let $p\geq 1$. Then there  exists  a constant $c >0$ such that for all integer $n \geq 1$, 
\begin{equation}
 \mathbb{E}[\sup_{0\leq t \leq T}|X^n_t|^{p}] \leq  c. \label{eqbb2}
  \end{equation}
\end{lem}

\noindent Proof of Lemma \ref{lembb1}.
Recall that, for all  $t\leq T$.,
\begin{equation}\label{eq:1biis}
X_t^n  =  \Phi_t^n  +  \alpha M_t^n  +  \beta I_t^n.
\end{equation}
Standard arguments,  using Burkholder-Davis-Gundy's inequality for the stochastic integral, H\"{o}lder's inequality, the linear growth condition (\ref{estimate:1}) satisfied by $b$ and $\sigma$, one show that 

\begin{equation*}
\mathbb{E}\sup_{0\leq t \leq T}\left|\int_{0}^{t}\sigma\left(s, X^n_{\eta_n(s)}\right)dW_s\right|^{p} \leq c  \int_{0}^{T}(1 + \mathbb{E}\sup_{0\leq u \leq s}\left|X^n_{u}\right|^{p})ds 
\end{equation*}
and 
\begin{equation*}
\mathbb{E}\sup_{0\leq t \leq T}\left|\int_{0}^{t}b\left(s, X^n_{\eta_n(s)}\right)ds\right|^{p} \leq c \int_{0}^{T}(1 + \mathbb{E}\sup_{0\leq u \leq s}\left|X^n_{u}\right|^{p})ds .
\end{equation*}
As a consequence we have the estimate
\begin{equation}\label{phi:est}
\mathbb{E}\sup_{0\leq t \leq T}|\Phi^n_t|^{p} \leq c \int_{0}^{T}(1 + \mathbb{E}\sup_{0\leq u \leq s}\left|X^n_{u}\right|^{p})ds .
\end{equation}

\noindent For  the term $M^n_t$, thanks to equation (\ref{eqbb7}) we have 
\begin{eqnarray}\label{max:est}
(1 - \alpha)|M^n_t|  & \leq & \max_{0 \leq s \leq t}|\Phi^n_t|  + \frac{|\beta|}{1 - \beta}\left(\max_{0 \leq s \leq \eta^{+}_n(t)}|\Phi^n_s|  + |\alpha|\max_{0 \leq s \leq \eta^{+}_n(t)}|M^n_{\eta^{+}_n(s)}|\right) \nonumber\\
& \leq & \left(1+\frac{|\beta|}{1 - \beta}\right) \max_{0 \leq s \leq t}|\Phi^n_s|  + \frac{|\alpha\beta|}{1-\beta}  \max_{0 \leq s \leq t}|M^n_{s}|.\nonumber
\end{eqnarray}
Therefore
\begin{equation}\label{min:est1}
(1 - |\rho|)\max_{0 \leq s \leq t}|M^n_s| \leq \frac{1}{1-\alpha}\left(1+\frac{|\beta|}{1 - \beta}\right)\max_{0 \leq s \leq t}|\Phi^n_s|, 
\end{equation}
\noindent Concerning  the term $I^n_t$, arguing as before and using this time the expression (\ref{eqbb7bis}) instead of (\ref{eqbb7}), we obtain
\begin{equation}\label{min:est}
(1 - |\rho|)\max_{0 \leq s \leq t}|I^n_s| \leq \frac{1}{1-\alpha}\left(1+\frac{|\beta|}{1 - \beta}\right)\max_{0 \leq s \leq t}|\Phi^n_s|, 
\end{equation}
Combining (\ref{max:est}) and (\ref{min:est}) together with (\ref{phi:est}), and using convexity of the function $x\rightarrow x^p$ inequality we get
\begin{equation*}
\mathbb{E}\sup_{0\leq s \leq T}\left|X^n_{s}\right|^{p} \leq c\int_{0}^{T}\left(1 + \mathbb{E}\sup_{0\leq u \leq s}\left|X^n_{u}\right|^{p}\right)ds. 
\end{equation*}
Now, by using  Gronwall's lemma, we obtain (\ref{eqbb2}). This finishes the proof of Lemma \ref{lembb1}.
\begin{lem}\label{lembb2}
Let $ p>1$,  $\gamma \in (0, \frac{1}{2} - \frac{1}{p})$ and $0 \leq s \leq t \leq T$. Then, 

\begin{equation} \label{eqbb5}
\mathbb{E}\left[\sup_{0 \leq s' \leq t' \leq T \atop |s' - t'| \leq |t-s| }|\Phi_{t'}^n - \Phi_{s'}^n|^{p}\right] \leq c|t-s|^{p \gamma}.
\end{equation}
\end{lem}

\noindent Proof of Lemma \ref{lembb2}.

\noindent Using the inequality $(a + b)^{p} \leq c (a^p +   b^p)$, valid for  all $a$, $b \geq 0$, we write
\begin{eqnarray}\label{eqbb34}
&& \mathbb{E}\left[\sup_{0 \leq s' \leq t' \leq T \atop |s' - t'| \leq |t-s| }|\Phi_{t'}^n - \Phi_{s'}^n|^p \right] \nonumber\\
&& \leq  c\left\{\mathbb{E}\sup_{0 \leq s' \leq t' \leq T \atop |s' - t'| \leq |t-s| }\left|\int_{s'}^{t'}\sigma\left(u, X^n_{\eta_n(u)}\right)dW_u \right|^p + \mathbb{E}\sup_{0 \leq s' \leq t' \leq T \atop |s' - t'| \leq |t-s| }\left|\int_{s'}^{t'}b\left(u, X^n_{\eta_n(u)}\right)du\right|^{p}\right\}.
\end{eqnarray}
\noindent Let us estimate the two terms appearing in the right hand side of (\ref{eqbb34}). 

\vspace{0.2cm}
\noindent For the first term, we will use Komogorov-Centsov criterion, see for instance \cite{Revuz-Yor-2005}. Indeed, set
 $\displaystyle  \left\{Y_{t}:=\int_{0}^{t'}\sigma\left(u,X^n_{\eta_n(u)}\right)dW_u, \ 0 \leq t' \leq T \right\}$, and observe that by Burkholder-Davis-Gundy's inequality, H\"{o}lder's inequality, the linear growth majoration satisfied by $\sigma$ as well as the estimation (\ref{eqbb2}), we have for $ 0\leq s'\leq t' \leq T$,
\begin{eqnarray}
\mathbb{E}\left[|Y_{t'}  -  Y_{s'}|^{p}\right]  & \leq&  c_{p} \mathbb{E}\left(\int_{s'}^{t'}\sigma^2\left(u, X^n_{\eta_n(u)}\right)du\right)^{p/2}\leq  2 c_{p} K^{2}\mathbb{E}\left( 1 + \sup_{u}|X^n_{u}|^2\right)^{p/2} |t' - s'|^{p/2}\nonumber \\ &\leq & c  |t' - s'|^{p/2}.
\end{eqnarray}

\noindent Then, according to (\cite{Revuz-Yor-2005}, Theorem 2.1; see also \cite{BarYor} for refinements and extensions based on Garsia-Rodemuch-Rumsey lemma), we deduce that
 $$
 \mathbb{E}\left[\left(\sup_{0\leq s' \leq t' \leq T}\frac{|Y_{t'}  -  Y_{s'}|}{|t' -s'|^{\gamma}} \right)^{p} \right] < \infty, 
 $$
for all $0 < \gamma < \frac{1}{p}\left(\frac{p}{2} - 1\right) = \frac{1}{2}  - \frac{1}{p}$, and therefore
 \begin{eqnarray}\label{eqbb35}
 \mathbb{E}\left[\sup_{0 \leq s' \leq t' \leq T \atop |s' - t'| \leq |t-s| } |Y_{t'} - Y_{s'}|^p\right] \leq |t  - s|^{p \gamma}. 
 \end{eqnarray}
\noindent For the remaining term in (\ref{eqbb34}), H\"{o}lder's inequality, linear growth majoration satisfied by $b$ as well as the estimation (\ref{eqbb2}) yield

\begin{eqnarray}\label{eqbb36}
\mathbb{E}\sup_{0 \leq s' \leq t' \leq T \atop |s' - t'| \leq |t-s| }\left|\int_{s'}^{t'}b\left(u, X^n_{\eta_n(u)}\right)du\right|^{p} & \leq &  \mathbb{E}\sup_{0 \leq s' \leq t' \leq T \atop |s' - t'| \leq |t-s|}\int_{s'}^{t'}\left|b\left(u, X^n_{\eta_n(u)}\right)\right|^{p}du \times |t'-s'|^{p-1}\nonumber\\
& \leq &  c \left(1 +  \mathbb{E}\sup_{0\leq u \leq T}|X^n_u|^{p}\right)\times |t-s|^{p}\nonumber \\
& \leq & c\left|t-s\right|^{p}.
\end{eqnarray}

\noindent The conclusion of Lemma \ref{lembb2} follows then from (\ref{eqbb35}) and (\ref{eqbb36}). \hfill $\square$ 
\begin{rk} The ``natural'' way to estimate the first term appearing in the right hand-side of (\ref{eqbb34}), would be to apply BDG's inequality, and obtain the exponent $p/2$ rather than $p\gamma$ in the inequality  (\ref{eqbb5}). This strategy failed, however, due to the facts that the supremum is taken over the two parameters $s^{\prime}$ and $t^{\prime}$, and that the doubly process $(\int_{s^{\prime}}^{t^{\prime}} \sigma\left(u, X^n_{\eta_n(u)}\right)dW_u))_{s^{\prime}, t^{\prime}}$ is no more a martingale.  It is interesting to note that, we found it more challenging and difficult to obtain the bound $O(|t-s|^{p/2})$ instead of $O(|t-s|^{p\gamma})$. 
\end{rk}

\vspace{0.2cm}
\begin{lem} \label{lem:2-5}Let $p>2$, and $\gamma \in (0, \frac{1}{2}-\frac{1}{p})$. Then, there exists a constant $c>0$ such that for all  integer $n\geq 1$, and $s, t \in [0, T]$, we have
\begin{equation}\label{eqbb6}
\mathbb{E}\sup_{0 \leq s' \leq t' \leq T \atop |s' - t'| \leq |t-s| }|M_{t'}^n - M_{s'}^n|^{p} \leq c\left|t-s \right|^{p\gamma}\,\,\,\, \mbox{and}\,\,\, \,  \mathbb{E}\sup_{0 \leq s' \leq t' \leq T \atop |s' - t'| \leq |t-s| }|I_{t'}^n - I_{s'}^n|^{p} \leq c\left|t-s \right|^{p\gamma}.
\end{equation}
Furthermore, \begin{equation} \label{eqbb18}
\mathbb{E}\sup_{0 \leq s' \leq t' \leq T \atop |s' - t'| \leq |t-s|}|X_{t'}^n - X_{s'}^n|^{p} \leq c\left|t-s \right|^{p\gamma}.
\end{equation}
\end{lem}

\noindent Proof of Lemma \ref{lem:2-5}. 
 In view of (\ref{eqbb7}), we can write for  $ 0\leq s\leq t \leq T$,
\begin{eqnarray}
0 \leq M_t^n - M_s^n & = &  \frac{1}{1 - \alpha}\max_{0 \leq u \leq t}\left(\Phi^n_{u} +  \frac{\beta}{\beta - 1}\max_{0\leq v \leq \eta^{+}_{n}{(u)}}\left(-\Phi^n_{v}-\alpha M^n_{\eta^{+}_n(v)}\right)^{+}\right)^{+}\nonumber \\
&  & - \frac{1}{1 - \alpha}\max_{0 \leq u \leq s}\left(\Phi^n_{u} +  \frac{\beta}{\beta - 1}\max_{0\leq v \leq \eta^{+}_{n}{(u)}}\left(-\Phi^n_{v}-\alpha M^n_{\eta^{+}_n(v)}\right)^{+}\right)^{+}.\nonumber 
\end{eqnarray}
If the maximum on the interval $[0,t]$ is attained on $[0,s]$, the above difference is zero. Otherwise, there exists $r\in (s, t]$ such that 
$$(1 - \alpha)M_t^n = \left(\Phi^n_{r} +  \frac{\beta}{\beta - 1}\max_{0\leq v \leq \eta_n^{+}(r)}\left(-\Phi^n_{v}-\alpha M^n_{\eta_n^{+}(v)}\right)^{+}\right)^{+}.$$
\noindent Then 
\begin{eqnarray}
0\leq M_t^n - M_s^n & \leq &   \frac{1}{1 - \alpha}\left(\Phi^n_{r} +  \frac{\beta}{\beta - 1}\max_{0\leq v \leq \eta_n^{+}(r)}\left(-\Phi^n_{v}-\alpha M^n_{\eta_n^{+}(v)}\right)^{+}\right)^{+}\nonumber \\
&  & - \frac{1}{1 - \alpha}\left(\Phi^n_{s} +  \frac{\beta}{\beta - 1}\max_{0\leq v \leq \eta_n^{+}(s)}\left(-\Phi^n_{v}-\alpha M^n_{\eta_n^{+}(v)}\right)^{+}\right)^{+}\nonumber\\
& \leq &  \frac{1}{1 - \alpha}\left|\Phi^n_{r} - \Phi^n_{s}\right|  + \frac{|\beta|}{(1-\alpha)(1-\beta)}\left|\Phi^n_{r'} - \Phi^n_{\eta_n^{+}(s)} \right|\nonumber\\
&&+ |\rho|\left|M^n_{\eta_n^{+}(r')} - M^n_{\eta_n^{+}(\eta_n^{+}(s))} \right|, \label{ineq:max}
\end{eqnarray}
where $r'\in ]\eta^{+}_{n}(s) , \eta^{+}_{n}(r)]$ satisfies  $ \max_{0\leq v \leq \eta_n^{+}(r)}\left(-\Phi^n_{v}-\alpha M^n_{\eta_n^{+}(v)}\right)^{+} = (-\Phi^n_{r'}-\alpha M^n_{\eta_n^{+}(r')})^{+}$  if  $\eta^{+}_{n}(s)<\eta^{+}_{n}(r)$ and
$$ \max_{0\leq v \leq \eta_n^{+}(r)}\left(-\Phi^n_{v}-\alpha M^n_{\eta_n^{+}(v)}\right)^{+} > \max_{0\leq v \leq \eta_n^{+}(s)}\left(-\Phi^n_{v}-\alpha M^n_{\eta_n^{+}(v)}\right)^{+}.$$

\vspace{0.2cm}

 \noindent Because $r'\in ]\eta_n^{+}(s) , \eta_n^{+}(r)]$,  $r\in (s, t]$ and the function $x \rightarrow x^{+}$ is  Lipschitz continuous on $\mathbb{R}$, we have
 $$|\eta_n^{+}(s)-r'|:=\left|\left(s- \frac{1}{n}\right)^{+} - r'\right| \leq \left|\left(s- \frac{1}{n}\right)^{+} - \left(r- \frac{1}{n}\right)^{+} \right| \leq  r-s \leq t-s,$$   and then
 \begin{eqnarray}
 |\eta_n^{+}(r') - \eta_n^{+}(\eta_n^{+}(s))| & := &   \left|\left(r'- \frac{1}{n}\right)^{+} - \left(\left(s-\frac{1}{n}\right)^{+} -\frac{1}{n}\right)^{+}\right|\leq   \left|r'-\left(s-\frac{1}{n}\right)^{+}\right| \leq t-s. \nonumber
\end{eqnarray}
Therefore, taking into account of (\ref{ineq:max}), there exists a positive constant $c$ depending only on $\alpha$ and $\beta$ such that 
    
\begin{equation*}
 (1- |\rho|)\sup_{0 \leq s' \leq t' \leq T \atop |s' - t'| \leq |t-s|} |M_{t'}^n - M_{s'}^n|   \leq c \sup_{0 \leq s' \leq t' \leq T \atop |s' - t'| \leq |t-s|} |\Phi^n_{t'}  - \Phi^n_{s'}|.
\end{equation*}
Taking expectation and using (\ref{eqbb5}) yields
 \begin{equation*}
\mathbb{E}\sup_{0 \leq s' \leq t' \leq T \atop |s' - t'| \leq |t-s| } |M_{t'}^n - M_{s'}^n|^{p}   \leq c \mathbb{E}\sup_{0 \leq s' \leq t' \leq T \atop |s' - t'| \leq |t-s| } |\Phi^n_{t'}  - \Phi^n_{s'}|^{p} \leq c\left|t - s\right|^{p\gamma}.
\end{equation*}
This proves the first inequality in (\ref{eqbb6}).

The second inequality  follows by arguing in the same way as above and using (\ref{eqbb7bis}) instead of (\ref{eqbb7}). Finally, the inequality (\ref{eqbb18})  is obtained as an immediate consequence of (\ref{eqbb5}), (\ref{eqbb6})  together with  (\ref{eq:1biis}) and the  convexity of the function $x \rightarrow|x|^p$.   This finished the proof of Lemma \ref{lembb2}.
\begin{rk} Iterating the inequality (\ref{ineq:max}) and using the preceding discussion enable us to obtain, 
\begin{eqnarray*}
0\leq M_t^n - M_s^n & \leq &  \left(\frac{1}{1 - \alpha} + \frac{|\beta|}{(1-\alpha)(1-\beta)}\right)\left( \sum_{k\geq 0} |\rho|^{k}\right)\sup_{0 \leq s' \leq t' \leq T \atop |s' - t'| \leq |t-s|} |\Phi^n_{t'}  - \Phi^n_{s'}|\\
&= &\frac{1}{1-|\rho|} \left(\frac{1}{1 - \alpha} + \frac{|\beta|}{(1-\alpha)(1-\beta)}\right)\sup_{0 \leq s' \leq t' \leq T \atop |s' - t'| \leq |t-s|} |\Phi^n_{t'}  - \Phi^n_{s'}|.
\end{eqnarray*}
Whence another way to check the first inequality of Lemma \ref{lem:2-5} by taking supremum and expectation in this inequality .
\end{rk}
\vspace{0.2cm}

We now turn to the proof of Theorem \ref{thmbb1}.

\vspace{0.2cm}
 
\noindent \textbf{Proof of Theorem \ref{thmbb1}}.\\
 \textit{ Proof of the existence:}
 As already  mentioned, we will show that $X^n$ is a Cauchy sequence in the Banach space $\mathrm{L}^{p}(\Omega, C([0,T];\mathbb{R}))$. Let $m$ and $n$ be two integers with $m > n$,                                                                                                     and $t\in [0, T]$. According to  equation (\ref{eq:1biis}), we can write
\begin{eqnarray}\label{eqbb15}
X^n_t - X^m_t & = & (\Phi^n_t  - \Phi^m_t) + \alpha(M_t^{n} - M_t^{m})  +  \beta(I_t^{n} - I_t^{m}).
\end{eqnarray}  
Notice that again we will make use of the representations (\ref{eqbb7}) and (\ref{eqbb7bis}) to handle the terms $\alpha(M_t^{n} - M_t^{m})$ and $\beta(I_t^{n} - I_t^{m})$. The proof is performed into several steps.\\
\noindent \textbf{Step $1$.}  We first claim that
\begin{equation}\label{eqbb27}
E[\sup_{t\in [0, T]} |M_t^n - M^m_t|^{p}]    \leq  c E[\sup_{t\in [0, T]} \left|\Phi^n_{t} - \Phi^m_{t}\right|^{p}]  + \  c \left|\frac{1}{n} - \frac{1}{m}\right|^{p\gamma}. 
\end{equation}
 Indeed, by (\ref{eqbb7}), Lipschitz continuity of the function $x\rightarrow x^{+}$ and elementary properties of the maximum, it is obviously seen that
\begin{eqnarray}\label{eqbb09}
&&|M_t^n - M^m_t|  \leq \frac{1}{1 - \alpha}\max_{0\leq s \leq t}\left|\Phi^n_{s} - \Phi^m_{s}\right| \nonumber\\
&& +   \frac{|\beta|}{(1 - \alpha)(1 - \beta)}\max_{0\leq s \leq t}\left|\max_{0 \leq u \leq \eta_n^{+}(s)}(-\Phi^n_{u} - \alpha M^{n}_{\eta_n^{+}(u)})^{+} - \max_{0 \leq u \leq \eta_m^{+}(s)}(-\Phi^m_{u} - \alpha M^{m}_{\eta_m^{+}(u)})^{+} \right|.\nonumber\\
&& \label{est:m,n}
\end{eqnarray}
In order to estimate the last term of this inequality, we distinguish two cases:\\
\noindent $\bullet$ \underline{\textit{Case $1$}}: If

 \begin{eqnarray} \label{case:1}
& & \left|\max_{0 \leq u \leq \eta_n^{+}(s)}(-\Phi^n_{u} - \alpha M^{n}_{\eta_n^{+}(u)})^{+} - \max_{0 \leq u \leq \eta_m^{+}(s)}(-\Phi^m_{u} - \alpha M^{m}_{\eta_m^{+}(u)})^{+} \right|\nonumber\\
& = &  \max_{0 \leq u \leq \eta_n^{+}(s)}(-\Phi^n_{u} - \alpha M^{n}_{\eta_n^{+}(u)})^{+} - \max_{0 \leq u \leq \eta_m^{+}(s)}(-\Phi^m_{u} - \alpha M^{m}_{\eta_m^{+}(u)})^{+},
\end{eqnarray}
\noindent let $r\in [0, \eta_n^{+}(s)]$ be such that   $\max_{0 \leq u \leq \eta_n^{+}(s)}(-\Phi^n_{u} - \alpha M^{n}_{\eta_n^{+}(u)})^{+} = (-\Phi^n_{r} - \alpha M^{n}_{\eta_n^{+}(r)})^{+}$. \\ Since  
$ \eta_m^{+}(s) \geq \eta_n^{+}(s) \geq r \geq 0 $, we have  $$\max_{0 \leq u \leq \eta_m^{+}(s)}(-\Phi^m_{u} - \alpha M^{m}_{\eta_m^{+}(u)})^{+} \geq  (-\Phi^m_{r} - \alpha M^{m}_{\eta_m^{+}(r)})^{+},$$  and  then (\ref{case:1}) gives
\begin{eqnarray} \label{eqbb29}
&& \left|\max_{0 \leq u \leq \eta_n^{+}(s)}(-\Phi^n_{u} - \alpha M^{n}_{\eta_n^{+}(u)})^{+} - \max_{0 \leq u \leq \eta_m^{+}(s)}(-\Phi^m_{u} - \alpha M^{m}_{\eta_m^{+}(u)})^{+} \right|\nonumber\\
&& \leq  (-\Phi^n_{r} - \alpha M^{n}_{\eta_n^{+}(r)})^{+} - (-\Phi^m_{r} - \alpha M^{m}_{\eta_m^{+}(r)})^{+} \nonumber\\
&& \leq  |(-\Phi^n_{r} - \alpha M^{n}_{\eta_n^{+}(r)}) - (-\Phi^m_{r} - \alpha M^{m}_{\eta_m^{+}(r)})| \nonumber\\
&& \leq  |\Phi^m_{r}-\Phi^n_{r}| + |\alpha| \left| M^{m}_{\eta_m^{+}(r)} - M^{n}_{\eta_n^{+}(r)}\right| \nonumber\\
&&  \leq  |\Phi^m_{r}-\Phi^n_{r}| + |\alpha| \left| M^{m}_{\eta_m^{+}(r)} - M^{m}_{\eta_n^{+}(r)}\right|  + |\alpha| \left| M^{m}_{\eta_n^{+}(r)} - M^{n}_{\eta_n^{+}(r)}\right| \nonumber\\
&& \leq \sup_{0\leq u \leq t} |\Phi^m_{u}-\Phi^n_{u}| + c \sup_{0 \leq s' \leq t' \leq t \atop |s' - t'| \leq |\frac{1}{m}-\frac{1}{n}|} |\Phi^m_{t'}  - \Phi^m_{s'}| + |\alpha|\sup_{0\leq u \leq t}|M^{m}_u-M^{n}_u|,
\end{eqnarray}
where we have used (\ref{eqbb6}) in the last inequality.

\noindent $\bullet$ \underline{\textit{Case $2$}}: If
\begin{eqnarray}
& & \left|\max_{0 \leq u \leq \eta_n^{+}(s)}(-\Phi^n_{u} - \alpha M^{n}_{\eta_n^{+}(u)})^{+} - \max_{0 \leq u \leq \eta_m^{+}(s)}(-\Phi^m_{u} - \alpha M^{m}_{\eta_m^{+}(u)})^{+} \right|\nonumber\\
& = & \max_{0 \leq u \leq \eta_m^{+}(s)}(-\Phi^m_{u} - \alpha M^{m}_{\eta_m^{+}(u)})^{+} - \max_{0 \leq u \leq \eta_n^{+}(s)}(-\Phi^n_{u} - \alpha M^{n}_{\eta_n^{+}(u)})^{+},\nonumber
\end{eqnarray}

\noindent let $r\in [0, \eta_m^{+}(s)]$ be such that   $\max_{0 \leq u \leq \eta_m^{+}(s)}(-\Phi^m_{u} - \alpha M^{m}_{\eta_m^{+}(u)})^{+} = (-\Phi^m_{r} - \alpha M^{m}_{\eta_m^{+}(r)})^{+}$. If  $r\in [0, \eta_n^{+}(s)]$, we  proceed as in the Case $1$ to get (\ref{eqbb29}) . Otherwise, $r\in (\eta_n^{+}(s), \eta_m^{+}(s)]$ and  we have 
\begin{eqnarray}\label{eeeq:2}
& & \left|\max_{0 \leq u \leq \eta_n^{+}(s)}(-\Phi^n_{u} - \alpha M^{n}_{\eta_n^{+}(u)})^{+} - \max_{0 \leq u \leq \eta_m^{+}(s)}(-\Phi^m_{u} - \alpha M^{m}_{\eta_m^{+}(u)})^{+} \right|\nonumber\\
& \leq & (-\Phi^m_{r} - \alpha M^{m}_{\eta_m^{+}(r)})^{+} - (-\Phi^n_{\eta_n^{+}(s)} - \alpha M^{n}_{\eta_n^{+}(\eta_n^{+}(s))})^{+} \nonumber\\
& \leq & |(-\Phi^m_{r} - \alpha M^{m}_{\eta_m^{+}(r)})- (-\Phi^n_{\eta_n^{+}(s)} - \alpha M^{n}_{\eta_n^{+}(\eta_n^{+}(s))})| \nonumber\\
& \leq & |\Phi^n_{\eta_n^{+}(s)}-\Phi^m_{r}| + |\alpha| | M^{m}_{\eta_m^{+}(r)}- M^{n}_{\eta_n^{+}(\eta_n^{+}(s))}| \nonumber\\
& \leq & |\Phi^n_{\eta_n^{+}(s)}-\Phi^n_{r}| + |\Phi^n_{r}-\Phi^m_{r}| + |\alpha| | M^{m}_{\eta_m^{+}(r)}- M^{m}_{\eta_n^{+}(\eta_n^{+}(s))}| + |\alpha| | M^{m}_{\eta_n^{+}(\eta_n^{+}(s))}- M^{n}_{\eta_n^{+}(\eta_n^{+}(s))}| \nonumber\\
\end{eqnarray}
\noindent Notice that, since $0 \leq \eta_n^{+}(s) \leq r \leq \eta_m^{+}(s)$, we obtain 
$$|\eta_n^{+}(s) - r| \leq |\eta_m^{+}(s) - \eta_n^{+}(s)| \leq |\eta_m(s) - \eta_n(s)| \leq |\frac{1}{m} - \frac{1}{n}|,$$
and  
\begin{eqnarray}
& &|\eta_n^{+}(\eta_n^{+}(s)) - \eta_m^{+}(r)| \leq  |\eta_n(\eta_n^{+}(s)) - \eta_m(r)| = \left|\left(s-\frac{1}{n}\right)^{+} - \frac{1}{n} -(r -\frac{1}{m})\right|\nonumber\\
&& \leq \left|\left(s-\frac{1}{n}\right)^{+} - r\right| +  \left|\frac{1}{m} - \frac{1}{n}\right| \leq \left|\left(s-\frac{1}{m}\right)^{+} - \left(s-\frac{1}{n}\right)^{+}\right| +  \left|\frac{1}{m} - \frac{1}{n}\right| \leq 2 \left|\frac{1}{m} - \frac{1}{n}\right|.\nonumber
\end{eqnarray}
Therefore, by (\ref{eeeq:2}) and (\ref{eqbb6}), we deduce

\begin{spreadlines}{12pt}
\begin{alignat}{2} \label{eqbb28}
&  \left|\max_{0 \leq u \leq \eta_n^{+}(s)}(-\Phi^n_{u} - \alpha M^{n}_{\eta_n^{+}(u)})^{+} - \max_{0 \leq u \leq \eta_m^{+}(s)}(-\Phi^m_{u} - \alpha M^{m}_{\eta_m^{+}(u)})^{+} \right|\nonumber\\
& \leq  |\Phi^n_{\eta_n^{+}(s)}-\Phi^n_{r}| + |\Phi^n_{r}-\Phi^m_{r}| + |\alpha| | M^{m}_{\eta_m^{+}(r)}- M^{m}_{\eta_n^{+}(\eta_n^{+}(s))}| + |\alpha| | M^{m}_{\eta_n^{+}(\eta_n^{+}(s))}- M^{n}_{\eta_n^{+}(\eta_n^{+}(s))}| \nonumber\\ 
&\leq   \sup_{0 \leq s' \leq t' \leq t \atop |s' - t'| \leq |\frac{1}{m}-\frac{1}{n}|} |\Phi^n_{t'}  - \Phi^n_{s'}|  +  \sup_{0\leq u \leq t} |\Phi^n_{u}-\Phi^m_{u}|+ |\alpha|\sup_{0 \leq s' \leq t' \leq t \atop |s' - t'| \leq 2|\frac{1}{m}-\frac{1}{n}|} |M^m_{t'}  - M^m_{s'}|  + |\alpha| \sup_{0\leq u \leq t}|M^{m}_u-M^{n}_u|
\end{alignat}
\end{spreadlines}

\noindent Observe that all terms of (\ref{eqbb29}) appear in (\ref{eqbb28}). Therefore, by inserting this last inequality in (\ref{eqbb09}) and taking the supremum over all the interval $[0, T]$, we infer that 

\begin{spreadlines}{12pt}
\begin{alignat}{2} 
& (1 - |\rho|)\sup_{t\in [0, T]} |M_t^n - M^m_t|\nonumber\\
&\leq   \sup_{0 \leq s' \leq t' \leq T \atop |s' - t'| \leq |\frac{1}{m}-\frac{1}{n}|} |\Phi^n_{t'}  - \Phi^n_{s'}|  +  c \sup_{t\in [0, T]} |\Phi^n_{t}-\Phi^m_{t}|+ |\alpha|\sup_{0 \leq s' \leq t' \leq T \atop |s' - t'| \leq 2|\frac{1}{m}-\frac{1}{n}|} |M^m_{t'}  - M^m_{s'}|.\nonumber 
\end{alignat}
\end{spreadlines}

\noindent Exponentiating both sides of the above inequality with $p$, and taking the expectation yield 

\begin{spreadlines}{12pt}
\begin{alignat}{2} 
& (1 - |\rho|)\mathbb{E}\left(\sup_{t\in [0, T]}|M_t^n - M^m_t|^p\right)\nonumber\\
&\leq \mathbb{E}\left(\sup_{0 \leq s' \leq t' \leq T \atop |s' - t'| \leq |\frac{1}{m}-\frac{1}{n}|} |\Phi^n_{t'}  - \Phi^n_{s'}|^{p}\right)  +  c \mathbb{E}\left(\sup_{t\in [0, T]} |\Phi^n_{t}-\Phi^m_{t}|^{p} \right) +  |\alpha|\mathbb{E}\left(\sup_{0 \leq s' \leq t' \leq T \atop |s' - t'| \leq 2|\frac{1}{m}-\frac{1}{n}|} |M^m_{t'}  - M^m_{s'}|^p\right)\nonumber \\
&\leq  c \mathbb{E}\left(\sup_{t\in [0, T]} |\Phi^n_{t}-\Phi^m_{t}|^{p} \right) +  \ c \left|\frac{1}{n} - \frac{1}{m}\right|^{p \gamma}, \nonumber 
\end{alignat}
\end{spreadlines}
where the last inequality follows readily from (\ref{eqbb5}) and (\ref{eqbb6}). This  proves the claim (\ref{eqbb27}).

\vspace{0.5cm}

\noindent \textbf{Step $2$.} By using (\ref{eqbb7bis}) instead of (\ref{eqbb3}), interchanging the roles of $M^n$ and $I^n$ in Step $1$ and making use of the second inequality of (\ref{eqbb6}),  we similarly obtain
\begin{equation} \label{eqbb25}
E[\sup_{t\in [0, T]} |I_t^n - I^m_t|^{p}]    \leq  c E[\sup_{t\in [0, T]} \left|\Phi^n_{t} - \Phi^m_{t}\right|^{p}]  + \  c \left|\frac{1}{n} - \frac{1}{m}\right|^{p\gamma}. 
\end{equation}
\\
\noindent \textbf{Step $3$.}  Let us prove the following statement
\begin{equation}\label{eqbb22}
\mathbb{E}\sup_{t\in [0, T]}|X_t^n - X^m_t|^{p}  \leq  c \left|\frac{1}{n} - \frac{1}{m}\right|^{p\gamma}.
\end{equation}

By convexity, we have from (\ref{eqbb15})
\begin{eqnarray}\label{eqbb26}
|X^n_t - X^m_t|^p & \leq  &3^{p-1} (|\Phi^n_t  - \Phi^m_t|^p + |\alpha(M_t^{n} - M_t^{m})|^p  +  |\beta(I_t^{n} - I_t^{m})|^p).
\end{eqnarray}  
Then, from Step $1$  and Step $2$, we only need to estimate  $\mathbb{E}\left(\sup_{t\in [0, T]}|\Phi^n_t - \Phi^m_t|^{p} \right)$.\\ \noindent It is obviously seen that
\begin{eqnarray}\label{eqbb13}
\mathbb{E}\left(\sup_{t\in [0, T]}|\Phi^n_t - \Phi^m_t|^{p} \right) & \leq & c \mathbb{E}\sup_{t\in [0, T]}\left|\int_{0}^{t}[\sigma\left(s, X^n_{\eta_n^{+}(s)}\right) - \sigma\left(s, X^m_{\eta_m^{+}(s)}\right)]dW_s \right|^{p}\nonumber\\
 &  & +  c \mathbb{E}\sup_{t\in [0, T]}\left|\int_{0}^{t}[b\left(s, X^n_{\eta_n^{+}(s)}\right) - b\left(s, X^m_{\eta_m^{+}(s)}\right)]ds \right|^{p}\nonumber\\
  &  \leq  & c \mathbb{E}\int_{0}^{T}\left|X^n_{\eta_n^{+}(s)} - X^m_{\eta_m^{+}(s)}\right|^p ds \nonumber\\
 &  \leq &  c\mathbb{E}\int_{0}^{T}\left|X^n_{\eta_n^{+}(s)} - X^m_{\eta_n^{+}(s)}\right|^p ds  + c \mathbb{E}\int_{0}^{T}\left|X^m_{\eta_n^{+}(s)} - X^m_{\eta_m^{+}(s)}\right|^p ds \nonumber\\
 &  \leq &  c \int_{0}^{T}\mathbb{E}\sup_{0\leq u \leq s }\left|X^n_{u} - X^m_{u}\right|^p ds  + c \left|\frac{1}{n} - \frac{1}{m}\right|^{p\gamma}, 
\end{eqnarray} 
where we have used (\ref{eqbb18}) in the last inequality. Now, consider the supremum over $[0, T]$, then take the expectation in (\ref{eqbb26}), and making use of the results of Step $1$ and Step $2$, we obtain
   
\begin{eqnarray}
\mathbb{E}\left(\sup_{t\in [0, T]}|X^n_t - X^m_t|^{p} \right) \leq   c \int_{0}^{T}\mathbb{E}(\sup_{0\leq u \leq s }\left|X^n_{u} - X^m_{u}\right|^p) ds  + c \left|\frac{1}{n} - \frac{1}{m}\right|^{p\gamma}.
\end{eqnarray}
Then Gronwall's  lemma implies  (\ref{eqbb22}) as was to be proved.

Now from Step $3$, we see then that the sequence $(X^n)$ is a Cauchy sequence in the complete space $\mathrm{L}^{p}(\Omega, C([0,T];\mathbb{R}))$, from which follows the existence of a continuous limit process $X : = (X(t), t\in [0, T])$ satisfying the property  (\ref{rate:1}) by letting $m$ goes to infinity in (\ref{eqbb22}). On the other hand, for every $\varepsilon >0$, $$\sum_{n} \mathbb{P}(\sup_{t\in [0, T]}|X^n_t - X_t| > \varepsilon)\leq \varepsilon^{-p}\sum_{n} n^{-p\gamma}< \infty, $$
for some $p>4$ and  $\gamma \in (\frac{1}{p}, \frac{1}{2}-\frac{1}{p})$. Thus, using  Borel-Cantelli lemma, we conclude that $(X^{n})$ converges to $X$ almost surely uniformly on $[0, T]$ as $n$ goes to infinity.

\noindent In order to complete the proof of Theorem \ref{thmbb1}, it remains to prove that $X$ is the unique solution of (\ref{eqbb1}). According to (\ref{eqbb13}) and(\ref{eqbb22}), Step $1$ and Step $2$, $(\Phi^n)$, $(M^n)$ and $(I^n)$ are Cauchy  sequences  in $\mathrm{L}^{p}(\Omega, C([0,T];\mathbb{R}))$. Let us denote by $\Phi$, $M$ and $N$ their limits respectively.  Passing to the limit, as $n$ goes to infinity, in  (\ref{eq:Phin}), (\ref{eqbb1bis}), (\ref{eqbb3}) and (\ref{eqbb4}) it is clearly seen that, almost surely, for all $t \geq 0$,
$$\displaystyle  \Phi_t  = \int_{0}^{t}\sigma\left(s, X_s\right)dW_s + \int_{0}^{t}b\left(s, X_s\right)ds\,\, \, \, \mbox{and}\,\, \,\, X_t= \Phi_t + \alpha M_t + \beta I_t.$$
\[M_t : = \frac{1}{1 - \alpha}\max_{0 \leq s \leq t}\left(\Phi_s + \beta I_{s} \right)^{+}, \,\,\, \mbox{and }\,\,\, I_t : = \frac{1}{\beta -1}\max_{0 \leq s \leq t}\left(-\Phi_s - \alpha M_{s} \right)^{+},\]
Therefore
\begin{equation}\label{eqbb71}
M_t = \frac{1}{1 - \alpha}\max_{0 \leq s \leq t}\left\{\left(\Phi_{s} +\frac{\beta}{\beta - 1}\max_{0\leq u \leq s}\left(-\Phi_{u} -\alpha M_{u}\right)^{+}\right)^{+} \right\}, 
\end{equation}  
and 
\begin{equation}\label{eqbb71bis}
I_t= \frac{1}{\beta - 1}\max_{0 \leq s \leq t}\left\{\left(-\Phi_{s} -\frac{\alpha}{ 1-\alpha}\max_{0\leq u \leq s}\left(\Phi_{u} +\beta I_{u}\right)^{+}\right)^{+} \right\}.
\end{equation}  
Now, we wish to show that $$M_t=\max_{0\leq s \leq t}X_s=: M_t^{X} \, \,\,\,\, \,\, \mbox{and}\,\,\,\,\,\,\, I_t=\min_{0\leq s \leq t}X_s=: I_t^{X}, \,t \geq 0.$$
This  will prove that $X$ satisfies (\ref{eqbb1}).  First, let us deal with $M_t$  and $M_t^X$. 

 On one hand, we have 
\begin{eqnarray}
M_t & = &  \frac{1}{1 - \alpha}\max_{0 \leq s \leq t}\left(\Phi_s + \beta I_{s} \right)^{+} \nonumber\\
& = &  \frac{1}{1 - \alpha}\max_{0 \leq s \leq t}\left(X_s - \alpha M_s \right)^{+}  \geq   \frac{1}{1 - \alpha}\left(X_t - \alpha M_t \right), 
\end{eqnarray}
which implies  $M_t \geq X_t$, and therefore $M_t \geq M^X_t$.
 On the other hand, let $t'\in [0, t]$ be such that 
$$
M_t = M_{t'} =  \frac{1}{1 - \alpha}\left(\Phi_{t'} +\frac{\beta}{\beta - 1}\max_{0\leq u \leq t'}\left(-\Phi_{u} -\alpha M_{u}\right)^{+} \right)^{+}.
$$ 
We notice two cases. If   $\Phi_{t'} +\frac{\beta}{\beta - 1}\max_{0\leq u \leq t'}\left(-\Phi_{u} -\alpha M_{u}\right)^{+} \leq 0$, that is, $M_{t'}=0$, then$$M^{X}_t \geq 0=M_{t'} = M_{t},$$ since $X_0=0$. Otherwise, we write 
\begin{eqnarray}
M_t^X &  = & \max_{0\leq s \leq t} X_s   =  \max_{0\leq s \leq t} (\Phi_s + \alpha M_s + \beta I_s) \geq  \Phi_{t'} + \alpha M_{t'} + \beta I_{t'}\nonumber\\
& = &  \alpha M_{t'} + \Phi_{t'} +   \frac{\beta}{\beta -1}\max_{0 \leq s \leq t'}\left(-\Phi_s - \alpha M_{s} \right)^{+} \nonumber\\
& = & \alpha M_{t'}  + (1-\alpha)M_{t'} =  M_{t'} = M_{t}. 
\end{eqnarray}
\noindent Thus, $M_t = M_t^X$.\\
\noindent   Concerning $I_t$ and $I_t^X$, we will proceed similarly. In fact, we have 
\begin{eqnarray}
I_t & = & \frac{1}{\beta -1}\max_{0 \leq s \leq t}\left(-\Phi_s - \alpha M_{s} \right)^{+} = \frac{1}{\beta -1}\max_{0 \leq s \leq t}\left(\beta I_s - X_s \right)^{+} \nonumber \\
& \leq &  \frac{1}{\beta -1}\left(\beta I_t - X_t \right),
\end{eqnarray} 
which is equivalent to $I_t \leq X_t$. Let $0\leq s \leq t$.  Since  the function $t\rightarrow-I_t$ is increasing, we get
$-X_s \leq  -I_s \leq -I_{t}$ and then $\max_{0\leq s \leq t}(-X_s) \leq -I_{t}$ which is equivalent to $I_{t} \leq I^{X}_{t}$.

\noindent  To show  that   $I_t^X \leq I_t$, we choose $s_0 \in [0, t]$ so that $$ (\beta -1)I_t=\left(-\Phi_{s_0} -\frac{\alpha}{ 1-\alpha}\max_{0\leq u \leq s_0}\left(\Phi_{u} +\beta I_{u}\right)^{+}\right)^{+},$$ and observe that 
$$ I_{t}=I_{s_0}=\frac{1}{\beta -1}\left(-\Phi_{s_0} -\frac{\alpha}{ 1-\alpha}\max_{0\leq u \leq s_0}\left(\Phi_{u} +\beta I_{u}\right)^{+}\right)^{+}.$$ Indeed, $I_{s_0} \geq I_t$ since the function $I_{\cdot}$ is non-increasing. On the other hand,  because $(-\Phi_{s_0} -\frac{\alpha}{ 1-\alpha}\max_{0\leq u \leq s_0}\left(\Phi_{u} +\beta I_{u}\right)^{+})^{+} \leq \max_{0 \leq s \leq s_0}\left(-\Phi_{s} -\frac{\alpha}{ 1-\alpha}\max_{0\leq u \leq s}\left(\Phi_{u} +\beta I_{u}\right)^{+} \right)^{+}$ and $\beta<1$, we obtain $I_{s_0}\leq I_{t}$. Thus $I_{t}=I_{s_0}$ as claimed. \\
\noindent  We also notice two cases.  If $-\Phi_{s_0} -\frac{\alpha}{ 1-\alpha}\max_{0\leq u \leq s_0}\left(\Phi_{u} +\beta I_{u}\right)^{+} \leq 0$,  then $$I^{X}_t \leq X_0  \leq 0 =I_t. $$Otherwise, we write   
\begin{eqnarray}
I_t^X & = & \min_{0 \leq s \leq t} X_s = - \max_{0 \leq s \leq t} (- \Phi_s - \alpha M_s - \beta I_s)\nonumber \\
 & = & -\max_{0 \leq s \leq t} \left(- \Phi_s -  \frac{\alpha}{1 - \alpha}\max_{0 \leq u \leq s}\left(\Phi_u + \beta I_{u} \right)^{+}  - \beta I_s \right)\nonumber \\
& \leq  & -\left(- \Phi_{s_0} -  \frac{\alpha}{1 - \alpha}\max_{0 \leq u \leq s_0}\left(\Phi_u + \beta I_{u} \right)  - \beta I_{s_0} \right)\nonumber \\
& = & -\left((\beta -1)I_{s_0}  - \beta I_{s_0} \right)= I_{s_0}=I_{t}.\nonumber\\
\end{eqnarray}
Summarizing, we have thus proved  $M_t=M^{X}_t$ and $I_t =I^{X}_t$ and therefore $$X_t= \Phi_t + \alpha M^{X}_t +\beta I^{X}_t.$$ This completes the proof of the existence part of Theorem \ref{thmbb1}.\\

\noindent \textit{ Proof of the uniqueness:} Now suppose that $X$ and $Y$ are two solutions to the SDE (\ref{eqbb1}) with the same driving Brownian motion $W$.  Then, 
\begin{equation}\label{eq1:uniq}
|X_t- Y_t|^{p}\leq 3^{p-1} ( |\Phi^{X}_{t}- \Phi^{Y}_{t}|^p+ |\alpha|^p |M^{X}_t-M^{Y}_t|^p + |\beta|^p |I^{X}_t -I^{Y}_t|^p,
\end{equation}
with $ \displaystyle  \Phi^{\omega}_t  = \int_{0}^{t}\sigma\left(s, \omega_s\right)dW_s + \int_{0}^{t}b\left(s, \omega_s\right)ds$,  where $\omega$ stand for $X$ or $Y$.
On one hand, using Skorohod's lemma,  we easily deduce:
\begin{equation}\label{eqbb111}
M^{\omega}_t = \frac{1}{1 - \alpha}\max_{0 \leq s \leq t}\left(\Phi^{\omega}_{s} +\frac{\beta}{\beta - 1}\max_{0\leq u \leq s}\left(-\Phi^{\omega}_{u} -\alpha M^{\omega}_{u}\right)^{+} \right)^{+}, 
\end{equation}  
and 
\begin{equation}\label{eqbb111bis}
I^{\omega}_t= \frac{1}{\beta - 1}\max_{0 \leq s \leq t}\left(-\Phi^{\omega}_{s} -\frac{\alpha}{ 1-\alpha}\max_{0\leq u \leq s}\left(\Phi^{\omega}_{u} +\beta I^{\omega}_{u}\right)^{+} \right)^{+}.
\end{equation}  
 On the other hand,  arguing as in Step $1$ one can  show that
\begin{equation}\label{eqbb1111bis}
(1 - |\rho|)\max_{t\in [0, t]} |M_s^X - M^Y_s| \leq c \max_{s\in [0, t]} |\Phi^{X}_s- \Phi^{Y}_s|,
\end{equation} 
and  
\begin{equation}\label{eqbb112bis} (1 - |\rho|)\max_{t\in [0, t]} |I_s^X - I^Y_s| \leq c \max_{s\in [0, t]} |\Phi^{X}_s- \Phi^{Y}_s|. 
\end{equation}
Moreover, by BDG's inequality and the Lipschitz conditions we get
\begin{equation}\label{eqbb113bis}
\mathbb{E}[\max_{s\in [0, t]} |\Phi^{X}_s- \Phi^{Y}_s|^{p}] \leq c \int_{0}^t \mathbb{E}(|X_s- Y_s|^p) ds.
\end{equation}
Taking expectation in (\ref{eq1:uniq}) and making use of (\ref{eqbb1111bis}), (\ref{eqbb112bis}) and (\ref{eqbb113bis}), we get
\[\mathbb{E}(|X_t-Y_t|^p \leq c\int_{0}^t \mathbb{E}(|X_s-Y_s|^p ds. )\]
Hence, $\mathbb{E}(|X_t-Y_t|^p =0$ by Gronwall's inequality. Thus the solution is unique.
\begin{rk} We  would like to emphasize that the choice of the initial condition of the solution $X$ in (\ref{eqbb1}) to be zero is important. Indeed,  the representations (\ref{eqbb111}) and (\ref{eqbb111bis}) used in the proof of the uniqueness part are only valid for $x_0=0$. Otherwise, the following approximation scheme could be used to show the existence result with the estimation (\ref{rate:1}) for any real random variable $x_0$ independent of $W$ such that $ \mathbb{E}(|x_0|^p)<\infty$:
\begin{equation*}
 X_t^n =  M_t^n = I_t^n = \frac{x_0}{1- \alpha - \beta},\,\,\,\, \text{for} \ -1 \leq t \leq 0,
\end{equation*}
and
\begin{equation*} 
 X_t^n  =  x_0 +  \int_{0}^{t}\sigma\left(s, X^n_{\eta_n(s)}\right)dW_s +  \int_{0}^{t}b\left(s, X^n_{\eta_n(s)}\right)ds + \alpha M_t^n  + \beta I_t^n, \,\,\text{for}\  t\in (0, T],
\end{equation*}
where $M_t^n$  and $I_t^n$ are respectively given by 
\begin{equation*}
M_t^n : = \frac{1}{1 - \alpha}\max_{0 \leq s \leq t}\left\{x_0+ \int_{0}^{s}\sigma\left(u, X^n_{\eta_n(u)}\right)dW_u +  \int_{0}^{s}b\left(u, X^n_{\eta_n(u)}\right)du + \beta I^n_{\eta^{+}_n(s)}\right\}, 
\end{equation*}
\begin{equation*}
I_t^n : = \frac{1}{\beta -1}\max_{0 \leq s \leq t}\left\{-x_0- \int_{0}^{s}\sigma\left(u, X^n_{\eta_n(u)}\right)dW_u -  \int_{0}^{s}b\left(u, X^n_{\eta_n(u)}\right)du -\alpha M^n_{\eta^{+}_n(s)}\right\}.
\end{equation*}   
\end{rk}
As a companion to Theorem \ref{thmbb1}, we have
\begin{prop}\label{prop:2} 
Let $1<p\leq 2$. Then, for every $\varepsilon \in (0, \frac{p}{2})$, there exists a constant $c>0$ such that 
\begin{equation} \label{rate:2}
\mathbb{E}\left(\sup_{t\in [0, T]}|X^n_t - X_t|^{p}\right) \leq \frac{c}{n^{p/2- \varepsilon}}.
\end{equation}
\end{prop}
\noindent Proof of Proposition \ref{prop:2}. Denote $a=\frac{p}{2}-\varepsilon$  and $\gamma=\frac{a}{p}$. Choose $p^{\prime}> \frac{2p}{p-2a}$ and observe that $p^{\prime}>2$ and $0<\gamma< \frac{1}{2}-\frac{1}{p^{\prime}}$, since $0<a<p/2$. Applying H\"{o}lder's inequality, we have 
\begin{eqnarray*}
\mathbb{E}\left[\sup_{t\in [0, T]}|X^n_t - X_t|^{p}\right] &= & \mathbb{E}\left[\left(\sup_{t\in [0, T]}|X^n_t - X_t|\right)^{p}\right]  \leq  \left[\mathbb{E}\left(\sup_{t\in [0, T]}|X^n_t - X_t|\right)^{p^{\prime}}\right]^{p/p^{\prime}}. 
\end{eqnarray*}
 In view of (\ref{rate:1}), we deduce
 \begin{eqnarray*}
\mathbb{E}\left[\sup_{t\in [0, T]}|X^n_t - X_t|^{p}\right] &\leq& c(n^{-p^{\prime} \gamma})^{p/p^{\prime}}=cn^{-p\gamma}=cn^{-a}.
\end{eqnarray*}
This finishes the proof of Proposition \ref{prop:2}.

As a consequence, we immediately obtain the following 
\begin{cor} Let $\varepsilon \in (0, 1)$, then there exists 
a constant $c>0$ such that 
\begin{equation*} 
\mathbb{E}\left(\sup_{t\in [0, T]}|X^n_t - X_t|^{2}\right) \leq \frac{c}{n^{1- \varepsilon}}.
\end{equation*}
\end{cor}

\section{A further extension: the non-Lipschitz case}
The main result of this Section is to show  that the scheme defined by (\ref{eqbb1bis}) still converges to the unique solution of (\ref{eqbb1}) with $x_0=0$, under the following  assumption:
\vspace{0.2cm}

\noindent \textbf{($\mathbf{H}2$)}: There exists $p>2$ such that for all $t \in [0, T]$ and $x, y \in \mathbb{R}$ we have:
\begin{equation*}
\label{no-lip} \left\{
\begin{array}{l}
|\sigma(t,x)-\sigma(t,y)|^p + |b(t,x)-b(t,y)|^p \leq  \rho(|x-y|^p),\,\,\, \\
\sup_{t\in [0, T]} (|\sigma(t, 0)| +|b(t, 0)|)< \infty, 
\end{array}
\right.
\end{equation*}
where $\rho$ is a continuous  non--decreasing and concave function $\rho(.):\R^+\to\R^+$ with $\rho(0)=0$, $\rho(u)>0$ for all $u>0$ and $\displaystyle\int_{0+}\frac{du}{\rho(u)}=+\infty$.\\
 Notice that assumption \textbf{($\mathbf{H2}$)} imply that  
\begin{equation}
\sup_{t\in [0, T]} |\sigma(t,x)|^{p}  \leq K ( 1 +  \rho\left(|x|^p\right))\, \, \, \mbox{and}\,\, \, \sup_{t\in [0, T]} |b(t,x)|^{p} \leq K ( 1 + \rho\left(|x|^p\right)).
\end{equation}
\vspace{0.2cm}

 \begin{thm}\label{thmbb2}  Assume that the functions $b$ and $\sigma$ are continuous and  satisfy  \textbf{($\mathbf{H}2$)}. Then, the SDE (\ref{eqbb1})  with $x_0=0$ has a unique strong solution. Moreover, for any $t>0$,  the sequence $(X^n)$ defined by (\ref{eqbb1bis}) converges in $\mathrm{L}^{p}(\Omega, C([0,T];\mathbb{R}))$. Furthermore,  
\begin{equation*} \label{limit}
\lim_{n \rightarrow +\infty}\mathbb{E}\left(\sup_{t\in [0, T]}|X^n_t - X_t|^{p}\right)=0.
\end{equation*}
\end{thm}
The proof of Theorem \ref{thmbb2} will require the following lemma.
\begin{lem}\label{lembb4}
Let $p> 2$. Under the assumption \textbf{$(\mathbf{H2})$}, there exists a positive constant $c$ such that

\begin{equation}
\sup_{n\geq 1}\mathbb{E}\left[ \sup_{t\in [0, T]} \left|X_t^n\right|^{p}\right] \leq c. 
\end{equation}

\end{lem}
\noindent Proof of Lemma (\ref{lembb4}). Recall (\ref{eq:1biis}) and write
\begin{equation}\label{eqbb33}
|X_t^n|^{p} \leq 3^{p-1} \left(|\Phi_t^n|^{p}  +  |\alpha|^{p} | M_t^n|^{p}  +  |\beta|^{p} |I_t^n|^{p}\right),
\end{equation}
for all $t\in[0, T]$. On one hand, we have 
\begin{eqnarray}\label{eqbb30}
&\displaystyle\mathbb{E}\sup_{0 \leq t \leq T }\left|\int_{0}^{t}\sigma\left(u, X^n_{\eta_n(u)}\right)dW_u \right|^p   \leq  c  \displaystyle\mathbb{E}\left(\int_{0}^{T}\sigma^{2}\left(u, X^n_{\eta_n(u)}\right)du\right)^{p/2}\nonumber\\
& \leq c\displaystyle  \mathbb{E}\int_{0}^{T}\left|\sigma\left(u, X^n_{\eta_n(u)}\right)\right|^{p}du \leq c \int_{0}^{T}\left(1 + \mathbb{E}\left[\rho\left(\left| X^n_{\eta_n(u)}\right|^p\right]\right)\right)du \nonumber\\
&  \leq  c \displaystyle \int_{0}^{T}\left(1 + \rho\left(\mathbb{E}\left| X^n_{\eta_n(u)}\right|^p\right)\right)du  \leq  c \int_{0}^{T}\left(1 + \rho\left(\mathbb{E}\sup_{0 \leq v \leq u }\left| X^n_{v}\right|^p\right)\right)du,
\end{eqnarray}
and similarly 
\begin{eqnarray}\label{eqbb31}
&\displaystyle\mathbb{E}\sup_{0 \leq t \leq T }\left|\int_{0}^{t}b\left(u, X^n_{\eta_n(u)}\right)du \right|^p \leq  c \int_{0}^{T}\mathbb{E}\left[\left|b\left(u, X^n_{\eta_n(u)}\right)\right|^p\right]du \nonumber\\
& \leq c\displaystyle  \int_{0}^{T}\left(1 + \mathbb{E}\left[\rho\left(\left| X^n_{\eta_n(u)}\right|^p\right]\right)\right)du
\leq c\int_{0}^{T}\left(1 + \rho\left(\mathbb{E}\sup_{0 \leq v \leq u }\left| X^n_{v}\right|^p\right)\right)du.
\end{eqnarray}
It follows from (\ref{eqbb30})  and (\ref{eqbb31}) that  
\begin{eqnarray}\label{eqbb32}
\mathbb{E}\left(\sup_{0 \leq t \leq T }\left|\Phi_t^n\right|^p\right) 
& \leq & c \int_{0}^{T}\left(1 + \rho \left(\mathbb{E}\sup_{0 \leq v \leq u }\left| X^n_{v}\right|^p\right)\right)du.
\end{eqnarray}

\noindent On the other hand,  (\ref{min:est1}) and (\ref{min:est})
 give 
$$ \displaystyle \mathbb{E}\left(\sup_{0 \leq t \leq T }\left|M_t^n\right|^p\right)  +  \mathbb{E}\left(\sup_{0 \leq t \leq T }\left|I_t^n\right|^p\right) \leq c\mathbb{E}\left(\sup_{0 \leq t \leq T }\left|\Phi_t^n\right|^p\right)\leq \displaystyle c \int_{0}^{T}\left(1 + \rho\left(\mathbb{E}\sup_{0 \leq v \leq u }\left| X^n_{v}\right|^p\right)\right)du 
.$$

\noindent Back to (\ref{eqbb33}), we therefore deduce that 
\begin{equation}
\mathbb{E}\left(\sup_{0 \leq t \leq T }\left|X_t^n\right|^p\right) \leq   c \int_{0}^{T}\left(1 + \rho\left(\mathbb{E}\sup_{0 \leq v \leq u }\left| X^n_{v}\right|^p\right)\right)du.
\end{equation}
Now, using  Bihari (also known as Bihari-LaSalle) inequality (see Appendix) we get $$\displaystyle  \sup_{n\geq 1}\mathbb{E}\left(\sup_{0 \leq t \leq T }\left|X_t^n\right|^p\right) \leq G^{-1}(G(cT) + c T) =c,$$
where $G(t):=\dint_{1}^{t} \frac{ds}{\rho(s)}$ and $G^{-1}$ denotes the inverse function of $G$.

 \noindent \textbf{Proof of Theorem \ref{thmbb2}}. \\
 We will again show that $(X^n)$ is a Cauchy sequence in the Banach space $\mathrm{L}^{p}(\Omega, C([0,T];\mathbb{R}))$. Let $m$ and $n$ be two integers with $m > n$  and $t\in [0, T]$. With the notations of the previous section, the results obtained in Step 1 and Step 2 together with (\ref{eqbb26}), we easily obtain
$$\mathbb{E}[\sup_{s\in [0, t]}|X^n_s - X^m_s|^{p} ] \leq c \mathbb{E}[\max_{s\in [0, t]} |\Phi^{n}_s- \Phi^{m}_s|^{p}]+  c \left|\frac{1}{n} - \frac{1}{m}\right|^{p\gamma},$$
for every $\gamma \in (0, \frac{1}{2}-\frac{1}{p})$.

 In order to  estimate $\mathbb{E}[\max_{s\in [0, t]} |\Phi^{n}_s- \Phi^{m}_s|^{p}]$, we use BDG's inequality, H\"{o}lder's inequality, assumption $(H2)$ and Jensen's inequality to get
\begin{eqnarray}\label{ineq:phi}
\mathbb{E}[\max_{s\in [0, t]} |\Phi^{n}_s- \Phi^{m}_s|^{p}] \leq c \displaystyle\mathbb{E}\left(\int_{0}^t \rho(|X_u^n-X_u^m|^p)du\right)\leq \displaystyle  \int_{0}^t \rho\left(\mathbb{E}[\sup_{u\in [0, s]}|X_u^n-X_u^m|^p]\right)ds.
\end{eqnarray}
Therefore
\begin{eqnarray*}
\mathbb{E}[\sup_{s\in [0, t]}|X^n_s - X^m_s|^{p} ]& \leq& \displaystyle  \int_{0}^t \rho\left(\mathbb{E}[\sup_{u\in [0, s]}|X^n_u - X^m_u|^{p} ] \right)ds +  c \left|\frac{1}{n} - \frac{1}{m}\right|^{p\gamma}.
\end{eqnarray*}
Notice that, in view of Lemma \ref{lembb4},  the continuous function $t \rightarrow \mathbb{E}[\sup_{s\in [0, t]}|X^n_s - X^m_s|^{p} ] $ is bounded. 
Then, applying  Bihari inequality, we obtain
\begin{eqnarray}\label{ineq:Bihari1}
\displaystyle\mathbb{E}[\sup_{s\in [0, t]}|X^n_s - X^m_s|^{p} ] \leq G^{-1}\left(G\left(\left|\frac{1}{n}-\frac{1}{m}\right|^{p\gamma}\right) + t\right).
\end{eqnarray}
From this and the fact that  $G\left(\left|\frac{1}{n}-\frac{1}{m}\right|^{p\gamma}\right) \rightarrow -\infty$ as $m, n \rightarrow +\infty$, since $\left|\frac{1}{n}-\frac{1}{m}\right|^{p\gamma} \rightarrow 0$ as well as $\dint_{0^{+}} \frac{ds}{\rho(s)}=\infty$
we get $$\lim_{m, n \rightarrow +\infty}\displaystyle\mathbb{E}[\sup_{s\in [0, t]}|X^n_s - X^m_s|^{p} ]=G^{-1}(-\infty)=0.$$
That is $(X^n)_n$ is a Cauchy sequence in the space $\mathrm{L}^{p}(\Omega, C([0,T];\mathbb{R}))$. Let $X$ be its limit. 

\noindent In order to complete the proof of Theorem \ref{thmbb2}, it remains then to prove that $X$ is the unique solution of the equation (\ref{eqbb1}). 

On one hand, combining (\ref{ineq:phi}) and (\ref{ineq:Bihari1}) and using the fact that  $\rho$ is non-decreasing yield
\[ \mathbb{E}[\max_{s\in [0, T]} |\Phi^{n}_s- \Phi^{m}_s|^{p}] \leq c T \rho\left(G^{-1}\left(G\left(\left|\frac{1}{n}-\frac{1}{m}\right|^{p\gamma}\right) + T\right)\right).
\]
 Thus $(\Phi^{n})_n$ is a Cauchy sequence in $\mathrm{L}^{p}(\Omega, C([0,T];\mathbb{R}))$, since $\rho(0)=0$.   On the other hand, by virtue of (\ref{eqbb27}) and  (\ref{eqbb25})  it follows that $(M^{n})_n$ and $(I^{n})_n$ are also Cauchy sequences in  $\mathrm{L}^{p}(\Omega, C([0,T];\mathbb{R}))$. Let $\Phi$ (resp.  $M$ and $I$) denotes the limit of $\Phi^{n}$, (resp. $M^n$ and $I^n$) in $\mathrm{L}^{p}(\Omega, C([0,T];\mathbb{R}))$.  So, by extracting a subsequence, that we continue to denote by $n$, we can assume without lost of generality  that  the sequences $X^n$, $\Phi^{n}$, $M^n$ and $I^n$ converge almost surely uniformly on $[0, T]$. Now, by letting $n$ tend to infinity in (\ref{eq:Phin}) and \eqref{eqbb1bis} and using  assumption $(H_2)$, the properties of $\rho$ and Bihari's inequality, one can easily prove that 
$$
\displaystyle  \Phi_t  = \int_{0}^{t}\sigma\left(s, X_s\right)dW_s + \int_{0}^{t}b\left(s, X_s\right)ds\,\, \, \, \mbox{and}\,\, \,\, X_t= \Phi_t + \alpha M_t + \beta I_t.$$

\noindent An analysis of the proof of Theorem  \ref{thmbb1} show that the arguments used after this step does not rely on Lipschitz properties of $b$ and $\sigma$  and  allow, mutatis matandis, to prove that $$M_t=\max_{0\leq s \leq t}X_s \, \,\,\,\, \,\, \mbox{and}\,\,\,\,\,\,\, I_t=\min_{0\leq s \leq t}X_s.$$
\noindent This concludes the proof of Theorem  \ref{thmbb2}.
\begin{examples}
Here we give some classical examples of functions $\rho$
satisfying  the above hypothesis $(H2)$ (see Yamada and Watanabe \cite{Yamada-Watanabe-1971}, and Barbu and Timisoara \cite{barbu-al-2000}). 
\begin{eqnarray*}
\hspace{-7.3cm}\rho_{1}(x) = x, \ \ \text{for} \ \ x \geq 0, 
\end{eqnarray*}
and
\begin{eqnarray*}
\rho_2(x) = \left\{
\begin{array}{l}
 0 \hspace{+1cm} if \quad    x=0 \\
 -x\log x    \hspace{+1cm} if \quad 0<x\leq\varepsilon, \ \ (\varepsilon \ \  \text{sufficiently small})\\
 \rho_1 (\varepsilon) +  \rho_1'(\varepsilon^{-})(x-\varepsilon)
 \hspace{0.5cm}if \quad   x>\varepsilon.
\end{array}
\right.
\end{eqnarray*}

\end{examples} 
\vspace{.5cm}
\section{Appendix}
\begin{lem}(Skorohod's lemma (\cite{Revuz-Yor-2005}, page 239))
Let $y$ be a real-valued continuous function on $[0, \infty[$ such that $y(0)\geq 0.$ Then there exists a unique pair $(z, k)$ of functions on $[0, \infty[$ such that 
\begin{itemize}
\item[i)]$ z= y + k$,
\item[ii)] $z$ is positive,
\item[iii)] $k$ is increasing, continuous, vanishing at zero, and is flat off $\{t\geq 0,, y(t)=0\}$; i.e. $\dint_{0}^{\infty} 1_{\{y(s)>0\}}dk(s)=0$. 
\end{itemize}
The function $k$ is give by 
\[ k(t)= \max_{0\leq s \leq t}((-y(s))^{+}), \, \,\, \mbox{for all}\, \,  t\geq 0.
\]
\end{lem}

\begin{lem}(Bihari's Lemma (\cite{Bihari1956}, pages 83-85), \cite{LaSalle49})
Let $T>0$ and $a \geq 0$. Assume that  $H$ and $G$ are two
nonnegative and measurable functions on $[0,T]$, and $\eta $ a
nonnegative continuous and nondecreasing function such that $\eta(r)
>0$ for $r>0$. If
\begin{displaymath}
H(t) \leq a + \int_{0}^{t} H(s)\eta(G(s))ds, \quad \mbox{for all}\
t\in[0,T],
\end{displaymath}
then
\begin{displaymath}
H(t) \leq G^{-1}\left(G(a) + \int_{0}^{t}H(s)ds \right)
\end{displaymath}
for all $t\in[0,T]$ such that $G(a) + \int_{0}^{t}H(s)ds\in
Dom(G^{-1})$, where $G(r) = \int_{1}^{r}\frac{ds}{\eta(s)}$ and
$G^{-1}$ is the inverse function  $G$.
\end{lem}

\end{document}